\documentclass[10pt]{amsart}
\usepackage{amssymb,latexsym}
\input epsf
\newtheorem{theorem}{Theorem}[section]
\newtheorem{proposition}[theorem]{Proposition}
\newtheorem{corollary}[theorem]{Corollary}

\newtheorem{remark}[theorem]{Remark}
\newtheorem{lemma}[theorem]{Lemma}
\newtheorem{example}[theorem]{Example}

\newcommand{\RR}{{\mathbb R}}

\newcommand{\cC}{{\mathcal C}}
\newcommand{\eE}{{\mathcal E}}
\newcommand{\fF}{{\mathcal F}}
\newcommand{\lL}{{\mathcal L}}

\newcommand{\NN}{{\mathbb N}}
\newcommand{\ZZ}{{\mathbb Z}}
\newcommand{\QQ}{{\mathbb Q}}
\newcommand{\kk}{{\bf k}}
\renewcommand{\to}{\rightarrow}

\newcommand{\conv}{{\rm conv}}
\newcommand{\sm}{{\setminus}}
\newcommand{\des}{{\rm des}}
\newcommand{\erh}{{\rm Ehrhart}}
\newcommand{\pull}{{\rm pull}}

\begin{document}
\title[Magic squares and a conjecture of Stanley]{Ehrhart polynomials, 
simplicial polytopes, magic squares and a conjecture of Stanley}

\author{Christos~A.~Athanasiadis}
\address{\hskip-\parindent Christos~A.~Athanasiadis\\
Department of Mathematics\\
University of Crete\\
71409 Heraklion, Crete, Greece}
\email{caa@math.uoc.gr}
\date{December 1, 2003}
%
\begin{abstract}
\centerline{\emph{Dedicated to Richard Stanley on the occasion of his 
sixtieth birthday}}

\vspace{0.1 in}
It is proved that a certain symmetric sequence $(h_0, h_1,\dots,h_d)$ of 
nonnegative integers arising in the enumeration of magic squares of given
size $n$ by row sums or, equivalently, in the generating function of the 
Ehrhart polynomial of the polytope of doubly stochastic $n \times n$ 
matrices, is equal to the $h$-vector of a simplicial polytope and hence 
that it satisfies the conditions of the $g$-theorem. The unimodality of 
$(h_0, h_1,\dots,h_d)$, which follows, was conjectured by Stanley (1983). 
Several generalizations are given.
\end{abstract}

\maketitle

\section{Introduction}
\label{intro}

A \emph{magic square} is a square matrix with nonnegative integer entries 
having all line sums equal to each other, where a line is a row or a column. 
Let $H_n (r)$ be the number of $n \times n$ magic squares with line sums equal 
to $r$. The problem to determine $H_n (r)$ appeared early in the twentieth 
century \cite{Ma}. Since then it has attracted considerable attention 
within areas such as combinatorics, combinatorial and computational 
commutative algebra, discrete and computational geometry, probability 
and statistics \cite{ADG, BP, DG, Eh, JvR, Sp, St1, St2, St4, St5, SS}. 
It was conjectured by Anand, Dumir and Gupta \cite{ADG} and proved by 
Ehrhart \cite{Eh} and Stanley \cite{St1} (see also 
\cite[Section I.5]{St4} and \cite[Section 4.6]{St5}) that for any fixed 
positive integer $n$, the quantity $H_n (r)$ is a polynomial in $r$ of 
degree $(n-1)^2$. More precisely, the following theorem holds. 
\begin{theorem} {\rm (Stanley~\cite{St1, St2})}
For any positive integer $n$ we have
\begin{equation}
\sum_{r \ge 0} \, H_n (r) \, t^r = \frac{h_0 + h_1 t + \cdots + h_d t^d}
{(1 - t)^{(n-1)^2 + 1}},
\label{mag0} 
\end{equation}
where $d = n^2 - 3n + 2$ and the $h_i$ are nonnegative integers satisfying 
$h_0 = 1$ and $h_i = h_{d-i}$ for all $i$.
\label{thm0}
\end{theorem}

It is the first conjecture stated in \cite{St4} (see Section I.1 there) that 
the integers $h_i$ appearing in (\ref{mag0}) satisfy further the inequalities
\begin{equation}
h_0 \le h_1 \le \cdots \le h_{\lfloor d/2 \rfloor}.
\label{mag1} 
\end{equation}
In this paper we prove this conjecture by showing that $(h_0, h_1,\dots,h_d)$  
is equal to the $h$-vector of a $d$-dimensional simplicial polytope. Such 
vectors are known to be symmetric and unimodal and are characterized by 
McMullen's $g$-theorem \cite{Mc}; see \cite{BL, St9} and Section 2. 

A few comments on Stanley's conjecture and the proof given in this paper are 
in order. It is known that the sequence $(h_0, h_1,\dots,h_d)$ of Theorem
\ref{thm0} is a \emph{Gorenstein sequence}, meaning it is the $h$-vector of 
a standard, graded, Gorenstein commutative ring; see for instance 
\cite[Section I.13]{St4}. It is an important open problem to characterize 
Gorenstein sequences; see \cite[Section II.6]{St4}. In this direction it 
was originally conjectured by Stanley \cite{St6} that a sequence $(h_0, 
h_1,\dots,h_d)$ is Gorenstein if and only if it satisfies the conditions 
of the $g$-theorem but a counterexample was later given in \cite{St7}. Our 
result (Corollary \ref{cor:magic}) and its generalization to the 
enumeration of magic labelings of regular bipartite graphs (Corollary 
\ref{cor:graphs}) give an instance in which Stanley's original conjecture 
turns out to be true. Another such instance, in which the entries of 
$(h_0, h_1,\dots,h_d)$ count linear extensions of a naturally labeled poset 
by the number of descents, was given recently by Reiner and Welker \cite{RW}. 
The polynomial $H_n (r)$ is the Ehrhart polynomial (see Section \ref{back})
of the Birkhoff polytope of doubly stochastic $n \times n$ matrices. We will 
state our main result in the context of Ehrhart polynomials of integer 
polytopes (Theorem \ref{cor:main}) as well as that of enumerating solutions 
to systems of linear homogeneous Diophantine equations (Corollary
\ref{thm:main}) and will show that both situations of Theorem \ref{thm0} 
and \cite{RW} appear as special cases. 

This paper was largely motivated by the work \cite{RW} of Reiner and Welker. 
I am grateful to Volkmar Welker for encouraging discussions and to Jes\'us
DeLoera, Victor Reiner, Francisco Santos and Richard Stanley for helpful 
suggestions.

\section{Background}
\label{back}

In this section we review some basic definitions and background on convex 
polytopes and their face numbers, triangulations and Ehrhart polynomials. 
We refer the reader to the texts by Stanley \cite{St4, St5}, Sturmfels 
\cite{Stu} and Ziegler \cite{Zi} for more information on these topics.   
We denote by $\NN$ the set of nonnegative integers.

\vspace{0.1 in}
\noindent
{\bf Face enumeration.} 
Given a finite (abstract or geometric) simplicial complex $\Delta$ of 
dimension 
$d-1$, let $f_i$ denote the number of $i$-dimensional faces of $\Delta$, so 
that $(f_0, f_1,\dots,f_{d-1})$ is the $f$-\emph{vector} of $\Delta$. The 
polynomial
\begin{equation}
\sum_{i=0}^d \, f_{i-1} (x-1)^{d-i} = \sum_{i=0}^d \, h_i x^{d-i}, 
\label{h-vec}
\end{equation}
where $f_{-1} = 1$ unless $\Delta$ is empty, is the $h$-\emph{polynomial} 
of $\Delta$, denoted $h(\Delta, t)$. The $h$-\emph{vector} of $\Delta$ is 
the sequence $(h_0, h_1,\dots,h_d)$ defined by (\ref{h-vec}). 

A \emph{polytopal complex} $\fF$ \cite[Section 8.1]{Zi} is a finite, 
nonempty collection of convex polytopes such that (i) any face of a polytope 
in $\fF$ is also in $\fF$ and (ii) the intersection of any two polytopes 
in $\fF$ is either empty or a face of both. The elements of $\fF$ are its 
\emph{faces} and those of dimension $0$ are its \emph{vertices}. The 
\emph{dimension} of $\fF$ is the maximum dimension of a face. The complex 
$\fF$ is \emph{pure} if all maximal faces of $\fF$ have the same dimension. 
The collection $\fF(P)$ of all faces of a polytope $P$ and the collection 
$\fF(\partial P)$ of its proper faces are pure polytopal complexes called 
the \emph{face complex} and \emph{boundary complex} of $P$, respectively. 
Thus $P$ is \emph{simplicial} if $\fF(\partial P)$ is a simplicial complex. 
The $h$-vectors of boundary complexes of simplicial polytopes are 
characterized by McMullen's $g$-\emph{theorem} \cite{Mc} \cite[Section 
III.1]{St4} \cite[Section 8.6]{Zi} as follows. A sequence $(g_0, 
g_1,\dots,g_\ell)$ of nonnegative integers is said to be an 
$M$-\emph{vector} if 
\begin{enumerate}
\itemsep=0pt
\item[(i)] $g_0 = 1$ and

\item[(ii)] $0 \le g_{i+1} \le g_i^{(i)}$ for $1 \le i \le \ell - 1$, 
\end{enumerate}
where $0^{(i)} = 0$ and
\[ n^{(i)} = {k_i + 1 \choose i+1} + {k_{i-1} + 1 \choose i} + \cdots + {k_j 
+ 1 \choose j+1} \]
for the unique representation
\[ n = {k_i \choose i} + {k_{i-1} \choose i-1} + \cdots + {k_j \choose j} \]
with $k_i > k_{i-1} > \cdots > k_j \ge j \ge 1$, if $n \ge 1$. A sequence 
$(h_0, h_1,\dots,h_d)$ of nonnegative integers is the $h$-vector of the 
boundary complex of a $d$-dimensional simplicial polytope if and only if
\begin{enumerate}
\itemsep=0pt
\item[(i)] $h_i = h_{d-i}$ for all $i$ and
\item[(ii)] $(h_0, h_1 - h_0,\dots,h_{\lfloor d/2 \rfloor} - h_{\lfloor d/2 
\rfloor - 1})$ is an $M$-vector.
\end{enumerate}
In particular $(h_0, h_1,\dots,h_d)$ is symmetric and unimodal and hence
satisfies the inequalities (\ref{mag1}), known as the \emph{Generalized 
Lower Bound Theorem} for simplicial polytopes.
 
\vspace{0.1 in}
\noindent
{\bf Triangulations and Ehrhart polynomials.} 
A \emph{triangulation} of a polytopal complex $\fF$ is a geometric 
simplicial complex $\Delta$ with vertices those of $\fF$ and underlying 
space equal to the union of the faces of $\fF$, such that every maximal
face of $\Delta$ is contained in a face of $\fF$. A triangulation 
of the face complex $\fF(P)$ of a polytope $P$ is simply called a 
triangulation of $P$. 

For any set $\sigma$ consisting of vertices of the polytopal complex $\fF$ 
we denote by $\fF \sm \sigma$ the subcomplex of faces of $\fF$ which do not 
contain any of the vertices in $\sigma$ and write $\fF \, \sm v$ for $\fF 
\sm \sigma$ if $\sigma$ consists of a single vertex $v$. Given a linear 
ordering $\tau = (v_1, v_2,\dots,v_p)$ of the set of vertices of $\fF$ we 
define the \emph{reverse lexicographic triangulation} or \emph{pulling 
triangulation} $\Delta (\fF) = \Delta_{\tau} (\fF)$ with respect to $\tau$ 
\cite{St3} \cite{Le} \cite[p.~67]{Stu} as $\Delta (\fF) = \{v\}$ if 
$\fF$ consists of a single vertex $v$ and
\[ \Delta (\fF) = \Delta (\fF \, \sm v_p) \, \cup \, \bigcup_F \ 
\{\conv(\{v_p\} \cup G): G \in \Delta (\fF(F)) \cup \{\emptyset\}\} \]
otherwise, where the union runs through the facets $F$ not containing 
$v_p$ of the maximal faces of $\fF$ which contain $v_p$ and $\Delta (\fF 
\, \sm v_p)$ and $\Delta (\fF(F))$ are defined with respect to the 
linear orderings of the vertices of $\fF \, \sm v_p$ and $F$, 
respectively, induced by $\tau$. Equivalently, for $i_0 < i_1 < \cdots 
< i_t$ the set $\{v_{i_0}, v_{i_1},\dots,v_{i_t}\}$ is the vertex set of 
a maximal simplex of $\Delta_{\tau} (\fF)$ if there exists 
a maximal flag $F_0 \subset F_1 \subset \cdots \subset F_t$ of faces of 
$\fF$ such that $v_{i_j}$ is the last vertex of $F_j$ with respect to 
$\tau$ for all $j$ and $v_{i_j}$ is not a vertex of $F_{j-1}$ for $j \ge 
1$. A different way to define $\Delta_{\tau} (\fF)$ is the following. For 
any vertex $v$ of $\fF$ let 
\[ \pull_v (\fF) = (\fF \, \sm v) \, \cup \, \bigcup_F \ \{\conv(\{v\} 
\cup G): G \in \fF(F) \cup \{\emptyset\}\}, \]
where the union runs through the facets $F$ not containing $v$ of the 
maximal faces of $\fF$ which contain $v$. If $\fF_0 = \fF$ and $\fF_i = 
\pull_{v_{p-i+1}} (\fF_{i-1})$ for $1 \le i \le p$ then $\fF_p$ is a 
triangulation of $\fF$ which coincides with $\Delta_{\tau} (\fF)$.
It follows from \cite[Theorem 2.5.23]{McS} (see also \cite[p.~80]{Gr}) that 
if $\fF$ is the boundary complex of a polytope $P$ then $\pull_v (\fF)$ is 
the boundary complex of another polytope, obtained from $P$ by moving its 
vertex $v$ beyond the hyperplanes supporting exactly those facets of $P$ 
which contain $v$. This observation implies the following lemma. 
\begin{lemma}
The reverse lexicographic triangulation of the boundary complex of a polytope
with respect to any ordering of its vertices is abstractly isomorphic to the 
boundary complex of a simplicial polytope of the same dimension.
\label{cor:rev}
\end{lemma}

A convex polytope $P \subseteq \RR^q$ is said to be a \emph{rational} or an 
\emph{integer} polytope if all its vertices have rational or integer 
coordinates, respectively. It is called a \emph{0-1 polytope} if all its 
vertices are 0-1 vectors in $\RR^q$. If $P$ is rational then the function 
defined for nonnegative integers $r$ by the formula
\[ \erh(P, r) = \# \, (r P \cap \ZZ^q) \]
is a quasi-polynomial in $r$, called the \emph{Ehrhart quasi-polynomial} of 
$P$ \cite[Section 4.6]{St5}. If $P$ is an integer polytope then this 
quasi-polynomial is actually a polynomial in $r$. Let $A \subseteq \RR^q$ be
the affine span of the integer polytope $P$. A triangulation $\Delta$ of $P$ 
is called \emph{unimodular} if the vertex set of any maximal simplex of 
$\Delta$ is a basis of the affine integer lattice $A \cap \ZZ^q$. We denote 
by $\Delta_\tau$ the reverse lexicographic triangulation of an arbitrary 
polytope $P$ with respect to the ordering $\tau$ of its vertices. Following 
\cite{St3} we call such an ordering of the vertices of an integer polytope 
$P$ \emph{compressed} if $\Delta_{\tau}$ is unimodular and call $P$ itself
\emph{compressed} if so is any linear ordering of its vertices. The following 
lemma holds for any unimodular triangulation of $P$, although we will not 
need this fact here. 
\begin{lemma} {\rm (\cite[Corollary 2.5]{St3})}
If $P$ is an $m$-dimensional integer polytope in $\RR^q$ and $\tau$ is a
compressed ordering of its vertices then
\[ \sum_{r \ge 0} \, \erh(P, r) \, t^r = \frac{h(\Delta_\tau, t)}
{(1 - t)^{m+1}}.\]
\label{lem:ehr}
\end{lemma}

If $P \subseteq \RR^m$ is an $m$-dimensional polytope and $V$ is any linear 
subspace of $\RR^m$ then the \emph{quotient polytope} $P / V \subseteq \RR^m 
/ \, V$ is the image of $P$ under the canonical surjection $\RR^m \to \RR^m 
/ \, V$. This is a convex polytope in $\RR^m / \, V$ linearly isomorphic to 
the image $\pi(P)$ of $P$ under any linear surjection $\pi: \RR^m \to 
\RR^{m - \dim V}$ with kernel $V$. Recall that the simplicial join $\Delta_1 
\ast \Delta_2$ of two abstract simplicial complexes $\Delta_1$ and $\Delta_2$ 
on disjoint vertex sets has faces the sets of the form $\sigma_1 \cup 
\sigma_2$, where $\sigma_1 \in \Delta_1$ and $\sigma_2 \in \Delta_2$ and that
$h(\Delta_1 \ast \Delta_2, t) = h(\Delta_1, t) \, h(\Delta_2, t)$. The 
following proposition is essentially Proposition 3.12 in \cite{RW}.
\begin{proposition}
Let $P$ be an $m$-dimensional polytope in $\RR^m$ having a triangulation 
abstractly isomorphic to $\sigma \ast \Delta$, where $\sigma$ is the vertex 
set of a simplex not contained in the boundary of $P$. Let $V$ be the linear 
subspace of $\RR^m$ parallel to the affine span of $\sigma$.

The boundary complex of the quotient polytope $P / V \subseteq \RR^m / 
\, V$ is abstractly isomorphic to $\fF(P) \sm \, \sigma$ and inherits a 
triangulation abstractly isomorphic to $\Delta$.
\label{prop:rw}
\end{proposition}

\vspace{0.1 in}
\noindent
{\bf Two compressed polytopes.} (a) A real $n \times n$ matrix is said to 
be \emph{doubly stochastic} if all its entries are nonnegative and all its
rows and columns sum to 1. The set $P$ of all real doubly stochastic $n 
\times n$ matrices is a convex polytope in $\RR^{n \times n}$ of dimension 
$(n-1)^2$, called the \emph{Birkhoff polytope} \cite[Example 0.12]{Zi}. It 
follows from the classical Birkhoff-von Neumann theorem that the vertices 
of $P$ are the $n \times n$ permutation matrices, so that $P$ is a 0-1 
polytope. The Birkhoff polytope was shown to be compressed by Stanley 
\cite[Example 2.4 (b)]{St3} (see also \cite[Corollary 14.9]{Stu}). 

\noindent
(b) Let $\Omega$ be a poset (short for partially ordered set) on the 
ground set $[m] := \{1, 2,\dots,m\}$. Recall that an (order) ideal of 
$\Omega$ is a subset $I \subseteq \Omega$ for which $i <_\Omega j$ and 
$j \in I$ imply 
that $i \in I$. Let $\Omega^0$ be the poset obtained from $\Omega$ by 
adjoining a minimum element $0$. The \emph{order polytope} \cite{St8} 
of $\Omega$, denoted $O(\Omega)$, is the intersection of the hyperplane 
$x_0 = 1$ in $\RR^{m+1}$ with the cone defined by the 
inequalities $x_i \ge x_j$ for $i < j$ in $\Omega^0$ and $x_i \ge 0$ for 
all $i$. Thus $O(\Omega)$ is an $m$-dimensional convex polytope. The 
vertices of $O(\Omega)$ are the characteristic vectors of the 
nonempty ideals of $\Omega^0$ \cite[Corollary 1.3]{St8} so, in particular, 
$O(\Omega)$ is a 0-1 polytope (see \cite[Theorem 1.2]{St8} for a complete 
description of the facial structure of $O(\Omega)$). Order polytopes were 
shown to be compresed by Ohsugi and Hibi \cite[Example 1.3 (b)]{OH}. 

\section{Special simplices}
\label{special}

Throughout this section $P$ denotes an $m$-dimensional convex polytope in 
$\RR^q$ with face complex $\fF(P)$. Let $\Sigma$ be a simplex spanned by $n$ 
vertices of $P$. We call $\Sigma$ a \emph{special simplex} in $P$ if each 
facet of $P$ contains exactly $n-1$ of the vertices of $\Sigma$. Note that, 
in particular, $\Sigma$ is not contained in the boundary of $P$. 
\begin{example} {\rm
Let $P$ be the polytope of real doubly stochastic $n \times n$ matrices. If 
$v_1, v_2,\dots,v_n$ are the $n \times n$ permutation matrices corresponding 
to the elements of the cyclic subgroup of the symmetric group generated by 
the cycle $(1 \ 2 \ \cdots \ n)$ (or any $n$ permutation matrices with pairwise 
disjoint supports) then $v_1, v_2,\dots,v_n$ are the vertices of a special 
simplex in $P$. Indeed, each facet of $P$ is defined by an equation of the 
form $x_{ij} = 0$ in $\RR^{n \times n}$ and misses exactly one of $v_1, 
v_2,\dots,v_n$. 
}
\label{ex1}
\end{example}
\begin{example} {\rm
Let $\Omega$ be a poset on the ground set $[m] := \{1, 2,\dots,m\}$ which is 
graded of rank $n-2$ (we refer to \cite[Chapter 3]{St5} for basic background 
and terminology on partially ordered sets) and $P = O(\Omega)$ be the order 
polytope of $\Omega$ in $\RR^{m+1}$. Let $\Omega^0$ be the poset obtained from 
$\Omega$ by adjoining a minimum element $0$ and for $1 \le i \le n$ let $v_i$ 
be the characteristic vector of the ideal of elements of $\Omega^0$ of rank 
at most $i-1$, so that $v_i$ is a vertex of $P$. Since a facet of $P$ is 
defined either by an equation of the form $x_i = x_j$ with $i < j$ in 
$\Omega^0$ and $i, j$ in successive ranks or by one of the form $x_i = 0$ for 
$i \in \Omega^0$ of rank $n-1$, it follows that $v_1, v_2,\dots,v_n$ are the 
vertices of a special simplex in $P$.  
}
\label{ex2}
\end{example}
\begin{lemma}
Suppose that $v_1, v_2,\dots,v_n$ are the vertices of a special simplex in 
$P$. If $F$ is a face of $P$ of codimension $k$ for some $1 \le k \le n-1$ 
and $F$ does not contain any of $v_1, v_2,\dots,v_k$ then $F$ must contain 
$v_i$ for all $k+1 \le i \le n$.
\label{lem1}
\end{lemma}
\begin{proof}
Let $\Sigma$ be the special simplex with vertices $v_1, v_2,\dots,v_n$. Any 
codimension $k$ face of a polytope can be written as the intersection
of $k$ facets, so we can write $F = F_1 \cap F_2 \cap \cdots \cap F_k$ where 
the $F_j$ are facets of $P$. For each $1 \le i \le k$ we have $v_i \notin F$ 
and hence $v_i \notin F_j$ for some $j = j_i$. Since $\Sigma$ is special the 
integers $j_1, j_2,\dots,j_k$ are all distinct and hence for each $1 \le j 
\le k$ we have $v_i \notin F_j$ for some $1 \le i \le k$, which in turn 
implies 
that $v_i \in F_j$ for all $k+1 \le i \le n$. It follows that $v_i \in F_1 
\cap F_2 \cap \cdots \cap F_k = F$ for all $k+1 \le i \le n$.  
\end{proof}
\begin{lemma}
Suppose that $\tau = (v_p, v_{p-1},\dots,v_1)$ is an ordering of the 
vertices of $P$ such that $\sigma = \{v_1, v_2,\dots,v_n\}$ is the vertex 
set of a special simplex in $P$. Let $\Delta$ be the abstract simplicial 
complex on $\{v_{n+1},\dots,v_p\}$ defined by the reverse lexicographic 
triangulation of $\fF(P) \sm \, \sigma$ with respect to $(v_p, 
v_{p-1},\dots,v_{n+1})$.  
\begin{enumerate}
\itemsep=0pt
\item[(i)] The reverse lexicographic triangulation $\Delta_{\tau}$ of $P$ 
is abstractly isomorphic to the simplicial join $\sigma \ast \Delta$.

\item[(ii)] $\Delta$ is abstractly isomorphic to the boundary complex of 
a simplicial polytope of dimension $m-n+1$.
\end{enumerate}   
\label{lem2}
\end{lemma}
\begin{proof}
(i) Let $\sigma_i = \{v_1,\dots,v_i\}$ for $0 \le i \le n$, so that $\sigma_0 
= \emptyset$ and $\sigma_n = \sigma$, and let $\Delta_i$ denote the abstract 
simplicial complex on the set $\{v_{i+1},\dots,v_p\}$ defined by the reverse 
lexicographic triangulation of $\fF(P) \sm \, \sigma_i$ with respect to the 
ordering $(v_p, v_{p-1},\dots,v_{i+1})$. To prove that $\Delta_0 = \sigma_n 
\ast \Delta_n$, which is the assertion in the lemma, we will prove that 
$\fF(P) \sm \, \sigma_i$ is pure $(m-i)$-dimensional and that $\Delta_0 = 
\sigma_i \ast \Delta_i$ for all $0 \le i \le n$ by induction on $i$. This 
is obvious for $i=0$ so let $1 \le i \le n$. By induction, any maximal face 
$F$ of $\fF(P) \sm \, \sigma_{i-1}$ is a codimension $i-1$ face of $P$. 
Since $F$ does not contain 
any of the vertices $v_1,\dots,v_{i-1}$, by Lemma \ref{lem1} we have $v_i \in 
F$. This implies that $\fF(P) \sm \, \sigma_i$ is pure $(m-i)$-dimensional and 
that $\Delta_{i-1} = v_i \ast \Delta_i$. The last equality and the induction 
hypothesis $\Delta_0 = \sigma_{i-1} \ast \Delta_{i-1}$ imply that $\Delta_0 = 
\sigma_i \ast \Delta_i$, which completes the induction. 

\noindent
(ii) Let $V$ be the linear subspace of $\RR^q$ parallel to the affine span 
of the vertices in $\sigma$ and $P / V$ be the corresponding quotient 
polytope of $P$, so that $P / V$ has dimension $m-n+1$. Part (i) and 
Proposition \ref{prop:rw} imply that $\Delta$ is abstractly isomorphic to a 
reverse lexicographic triangulation of the boundary complex of $P / V$. This 
is in turn isomorphic to the boundary complex of a simplicial polytope of 
dimension $m-n+1$ by Lemma \ref{cor:rev}.
\end{proof}

The following theorem is the key to the results in this paper. 
\begin{theorem}
Suppose that $P$ is an integer polytope and $\tau = (v_p, 
v_{p-1},\dots,v_1)$ is an ordering of its vertices such that: 
\begin{enumerate}
\itemsep=0pt
\item[(i)] $\tau$ is compressed and
\item[(ii)] $\{v_1, v_2,\dots,v_n\}$ is the vertex set of a special simplex 
in $P$.  
\end{enumerate}   
Then
\[ \sum_{r \ge 0} \, \erh(P, r) \, t^r = \frac{h(t)}{(1 - t)^{m+1}} \]
where $h(t) = h_0 + h_1 t + \cdots + h_d t^d$ is the $h$-polynomial of 
the boundary complex of a simplicial polytope $Q$ of dimension $d=m-n+1$, 
so that $h(t)$ satisfies the conditions in the $g$-theorem. 

In particular  $h_i = h_{d-i}$ for all $i$ and $1 = h_0 \le h_1 \le \cdots 
\le h_{\lfloor d/2 \rfloor}$.

Moreover, $Q$ can be chosen so that its boundary complex is abstractly 
isomorphic to the reverse lexicographic triangulation of $\fF(P) \sm \, 
\{v_1,\dots,v_n\}$ with respect to the ordering $(v_p, 
v_{p-1},\dots,v_{n+1})$.
\label{cor:main}
\end{theorem}
\begin{proof}
Let $\sigma = \{v_1, v_2,\dots,v_n\}$ and let $\Delta$ denote the reverse 
lexicographic triangulation of $\fF(P) \sm \, \sigma$ with respect to the 
ordering $(v_p, v_{p-1},\dots,v_{n+1})$. Lemma \ref{lem:ehr} guarantees 
that the proposed equation holds with $h(t) = h(\Delta_\tau, t)$. Part (i) 
of Lemma \ref{lem2} implies that 
$$h(\Delta_\tau, t) = h(\sigma \ast \Delta, t) = h(\sigma, t) \, h(\Delta, 
t) = h(\Delta, t),$$ since face complexes of simplices have $h$-polynomial 
equal to $1$, and the result follows from part (ii) of the same lemma. 
\end{proof}

We now apply Theorem \ref{cor:main} to the Birkhoff polytope and to order 
polytopes of graded posets. Observe that our theorem does not apply to all 
integer polytopes since 0-1 polytopes with no regular unimodular 
triangulations are known to exist \cite{OH2}.

\vspace{0.1 in}
\noindent
{\bf Magic squares and the Birkhoff polytope.} Let $P$ be the polytope of 
real doubly stochastic $n \times n$ matrices. Observe that the polynomial 
$\erh(P, r)$ coincides with the function $H_n (r)$ of Theorem \ref{thm0}. 
Since $P$ is a compressed integer polytope of dimension $(n-1)^2$, 
Theorem \ref{cor:main} and Example \ref{ex1} imply immediately the 
following corollary.
\begin{corollary}
For any positive integer $n$ we have
\[ \sum_{r \ge 0} \, H_n (r) \, t^r = \frac{h(t)} {(1 - t)^{(n-1)^2 + 1}} \]
where $h(t) = h_0 + h_1 t + \cdots + h_d t^d$ is the $h$-polynomial of the 
boundary complex of a simplicial polytope of dimension $d = n^2 - 3n + 2$, 
so that $h(t)$ satisfies the conditions in the $g$-theorem. 

In particular  $h_i = h_{d-i}$ for all $i$ and $1 = h_0 \le h_1 \le \cdots 
\le h_{\lfloor d/2 \rfloor}$.
\label{cor:magic}
\end{corollary}
In view of the last statement in Corollary \ref{thm:main}, the polytope in 
the previous corollary can be constructed by pulling in an arbitrary order 
the vertices of the quotient of $P$ with respect to the affine span of the 
vertices $v_1, v_2,\dots,v_n$, chosen explicitly as in Example \ref{ex1}. 

\vspace{0.1 in}
\noindent
{\bf Eulerian polynomials and equatorial spheres.} 
Let $\Omega$ be a graded poset on the ground set $[m] := \{1, 2,\dots,m\}$ 
of rank $n-2$. Let $\Omega_i$ be the set of elements of $\Omega$ of rank 
$i-1$ for $1 \le i \le n-1$ and $\lL(\Omega)$ be the set of \emph{linear 
extensions} of $\Omega$, meaning the set of permutations $w = (w_1, 
w_2,\dots,w_m)$ of $[m]$ for which $w_i <_\Omega w_j$ implies $i < j$. We 
assume that $\Omega$ is \emph{naturally labeled}, meaning that the identity 
permutation $(1, 2,\dots,m)$ is a linear extension. The 
$\Omega$-\emph{Eulerian polynomial} is defined as
\[ W(\Omega, t) = \sum_{w \in \lL(\Omega)} \ t^{\des (w)}\]
where 
\[ \des (w) = \# \, \{i \in [m-1]: w_i > w_{i+1} \} \]
is the number of descents of $w$. Following \cite{RW} we call a 
function $g: \Omega \to \RR$ \emph{equatorial} if $\min_{a \in \Omega} 
g(a) = 0$ and for each $2 \le i \le n-1$ there exist $a_{i-1} \in 
\Omega_{i-1}$ and $a_i \in \Omega_i$ such that $a_{i-1} <_{\Omega} 
a_i$ and $g(a_{i-1}) = g(a_i)$. An ideal $I$ or, more generally, a 
strictly increasing chain of ideals $I_1 \subset I_2 \subset \cdots 
\subset I_k$ in $\Omega$ is \emph{equatorial} if the characteristic 
function $\chi_I$ of $I$ or the sum $\chi_{I_1} + \chi_{I_2} + \cdots + 
\chi_{I_k}$, respectively, is equatorial. The \emph{equatorial complex} 
$\Delta_{eq} (\Omega)$, introduced in \cite{RW}, is the abstract simplicial 
complex on the vertex set of equatorial ideals of $\Omega$ whose simplices 
are the equatorial chains of ideals in $\Omega$. 

The following theorem is proved in Corollary 3.8 and Theorem 3.14 of \cite{RW}.
\begin{theorem} {\rm (Reiner--Welker~\cite{RW})}
Let $\Omega$ be a naturally labeled, graded poset on $[m]$ having $n-1$ ranks. 
The equatorial complex $\Delta_{eq} (\Omega)$ is abstractly isomorphic to the 
boundary complex of a simplicial polytope of dimension $d = m-n+1$ which has
$h$-polynomial equal to the $\Omega$-Eulerian polynomial $W(\Omega, t)$.

Hence $W(\Omega, t)$ satisfies the conditions in the $g$-theorem and, in 
particular, it has symmetric and unimodal coefficients.
\label{thm:eq}
\end{theorem}

Let $P$ be the order polytope of $\Omega$ and $\Omega^0$ be the poset 
obtained from $\Omega$ by adjoining a minimum element $0$. Recall that 
the vertices of $P$ are the characteristic vectors of the nonempty
ideals of $\Omega^0$. The order polytope comes with 
its canonical triangulation \cite{St8} \cite[Proposition 2.1]{RW}, which 
is a unimodular triangulation with maximal simplices bijecting to the 
linear extensions of $\Omega$. This canonical triangulation is in fact 
the reverse lexicographic triangulation of $O(\Omega)$ with respect to 
any ordering $(u_p, u_{p-1},\dots,u_1)$ of its vertices such that $i < j$ 
whenever the ideal of $\Omega^0$ defined by $u_i$ is strictly contained in 
that defined by $u_j$. We will use the following lemma.  
\begin{lemma}
Let $v_i$ be the characteristic vector of the ideal of elements of 
$\Omega^0$ of rank at most $i-1$ for $1 \le i \le n$. Let $\sigma = 
\{v_1,\dots,v_n\}$ and $\tau = (v_p,\dots,v_{n+1})$ be an ordering of 
the remaining vertices of $P$ such that $i < j$ whenever $i, j \ge n+1$ 
and the ideal defined by $v_i$ is strictly contained in that defined 
by $v_j$. 

The equatorial complex $\Delta_{eq} (\Omega)$ is the abstract simplicial 
complex defined by the reverse lexicographic triangulation of $\fF(P) \sm 
\, \sigma$ with respect to $\tau$. 
\label{lem:eq}
\end{lemma}
\begin{proof}
Let $\fF$ denote the face complex of $P$ and let $x_{\hat{1}} = 0$ by 
convention. The maximal faces of $\fF \sm \sigma$ are the faces of $P$ 
defined by systems of equations of the form $x_{i_s} = x_{j_s}$ for $0 \le s 
\le n-1$ where (i) $i_0 = 0$ and $j_{n-1} = \hat{1}$, (ii) $i_s \in \Omega_s$ 
for $1 \le s \le n-1$, $j_s \in \Omega_{s+1}$ for $0 \le s \le n-2$ and $i_s 
<_\Omega j_s$ for $1 \le s \le n-2$ and (iii) if $j_s = i_{s+1}$ for 
consecutive values $s = a, a+1,\dots,b-1$ of $s$ then the interval $[i_a, 
j_b]$ in $\hat{\Omega}$ consists only of the elements of the chain $i_a < 
i_{a+1} < \cdots < i_b < j_b$. The statement of the Lemma follows from the 
description of the maximal faces of a reverse lexicographic triangulation 
$\Delta(\fF)$ (see Section \ref{back}) and that of the maximal faces of 
$\Delta_{eq} (\Omega)$ (see \cite[Proposition 3.5]{RW}). We omit the details. 
\end{proof}

\vspace{0.1 in}
\noindent
\emph{Proof of Theorem \ref{thm:eq}.} Let $P$ be the order polytope of 
$\Omega$, as before. Observe that $\erh(P, r)$ is equal to the number of order 
reversing maps $\rho: \Omega \to \{0, 1,\dots,r\}$. It follows from 
\cite[Theorem 4.5.14]{St5} that
\begin{equation}
\sum_{r \ge 0} \, \erh(P, r) \, t^r = \frac{W(\Omega, t)} {(1 - t)^{m + 1}}.
\label{st-P}
\end{equation}
Let the vertices $v_1, v_2,\dots,v_p$ of $P$ and $\tau = (v_p,\dots,v_{n+1})$ 
be as in Lemma \ref{lem:eq}. We checked in Example \ref{ex2} that $v_1, 
v_2,\dots,v_n$ are the vertices of a special simplex in $P$. Since $P$ is a 
compressed integer polytope (see Section \ref{back}) Theorem \ref{cor:main} 
applies and we have
\[ \sum_{r \ge 0} \, \erh(P, r) \, t^r = \frac{h(t)} {(1 - t)^{m + 1}}, \]
where $h(t) = h_0 + h_1 t + \cdots + h_d t^d$ is the $h$-polynomial of a 
simplicial polytope of dimension $d=m-n+1$ having, in view of Lemma 
\ref{lem:eq}, boundary complex abstractly isomorphic to $\Delta_{eq} (\Omega)$.
Comparison with (\ref{st-P}) yields $h(t) = W(\Omega, t)$ and completes the 
proof.
\qed

\section{Rational polyhedral cones}
\label{linear}

In this section we state several corollaries of Theorem \ref{cor:main},
including a generalization of Corollary \ref{cor:magic} to magic labelings 
of bipartite graphs. 
Let $\Phi$ be a $p \times q$ integer matrix with rank $p$ and $E_\Phi$ be 
the monoid of vectors $x \in \NN^q$ satisfying the homogeneous system of 
linear equations $\Phi \, x = 0$. We assume that $E_\Phi$ is nonzero, so 
that $q \ge p + 1$, and set $m = q-p-1$. We denote by $\cC_\Phi$ the cone 
of vectors $x = (x_1, x_2,\dots,x_q) \in \RR^q$ satisfying $\Phi \, x = 0$ 
and $x_i \ge 0$ for all $i$. Let $L(x) = a_1 x_1 + a_2 x_2 + \cdots + a_q 
x_q$ be a linear functional on $\RR^q$ with $a_i \in \QQ$ for all $i$ such 
that the set 
\begin{equation}
P = \{x \in \cC_\Phi: L(x) = 1\}
\label{P}
\end{equation}
is nonempty and bounded. Thus $P$ is an $m$-dimensional convex polytope. 
Observe that  
\begin{equation}
\erh(P, r) = \# \, \{x \in E_\Phi: L(x) =  r\}
\label{f}
\end{equation}
for any nonnegative integer $r$. Let $\bar{E}_\Phi$ be the submonoid 
of $E_\Phi$ consisting of those elements having positive coordinates. 
We say that $\beta \in \bar{E}_\Phi$ is the \emph{unique minimal 
element} of $\bar{E}_\Phi$ if we have $\beta \le \gamma$ 
coordinatewise for all $\gamma \in \bar{E}_\Phi$.
\begin{corollary}
Let $P$ be an $m$-dimensional polytope, as in {\rm (\ref{P})}. If
\begin{enumerate}
\itemsep=0pt
\item[(i)] $P$ is an integer polytope, 

\item[(ii)] $\bar{E}_\Phi$ has a unique minimal element $\beta$ and

\item[(iii)] there exists a compressed ordering $\tau = (v_p, 
v_{p-1},\dots,v_1)$ of the vertices of $P$ such that $v_1 + v_2 + 
\cdots + v_n = \beta$ for some $n$
\end{enumerate}   
then the conclusion of Theorem \ref{cor:main} holds.
\label{thm:main}
\end{corollary}
\begin{proof}
In view of Theorem \ref{cor:main} it suffices to show that $\{v_1, 
v_2,\dots,v_n\}$ is the vertex set of a special simplex in $P$. Let 
$\beta = (\beta_1, \beta_2,\dots,\beta_q)$ and let $F$ be a facet of 
$P$, so that $F$ is defined by an equation of the form $x_k = 0$ for 
some $1 \le k \le q$. We need to show that exactly one of $v_1, 
v_2,\dots,v_n$ has positive $k$th coordinate. Clearly at least one 
of $v_1, v_2,\dots,v_n$ has this property, since $\beta_i > 0$ for
all $i$. Assume on the contrary that at least two $v_j$ have positive 
$k$th coordinate, say $v_1$ and $v_2$, so that $1 \le \gamma_k < 
\beta_k$ if $v_1 = (\gamma_1, \gamma_2,\dots,\gamma_q)$. Since $F$ 
is a facet of $P$ there exists a point $x = (x_1, x_2,\dots,x_q)$ in 
the affine span of $P$, which we may assume to be rational, 
satisfying $x_k < 0$ and $x_i > 0$ for all $i \neq k$. By replacing 
$x$ with a suitable integer multiple we find a point $\alpha = 
(\alpha_1, \alpha_2,\dots,\alpha_q) \in \ZZ^q$ satisfying $\Phi \, 
\alpha = 0$, $\alpha_k < 0$ and $\alpha_i > 0$ for all $i \neq k$. 
We may choose a nonnegative integer $t$ so that $0 < \alpha_k + 
\beta_k + t \gamma_k < \beta_k$ (with $t=0$ if $\alpha_k + \beta_k 
> 0$). Then $\alpha + \beta + t v_1$ is in $\bar{E}_\Phi$ and has $k$th 
coordinate strictly less than $\beta_k$, which contradicts the 
minimality of $\beta$.
\end{proof}

The condition that $\bar{E}_\Phi$ has a unique minimal element is
satisfied if $(1, 1,\dots,1) \in \RR^q$ is in $E_\Phi$ and is known 
to hold if and only if the semigroup ring $R_\Phi = \kk [E_\Phi]$, 
generated over a field $\kk$ by the monomials corresponding to 
elements of $E_\Phi$, is Gorenstein; see for instance \cite[Section 
I.13]{St4}. 
\begin{corollary}
Let $P$ be an $m$-dimensional polytope, as in {\rm (\ref{P})}. If
\begin{enumerate}
\itemsep=0pt
\item[(i)] $P$ is a compressed integer polytope, 
\item[(ii)] $R_\Phi$ is Gorenstein and
\item[(iii)] $E_\Phi$ is generated as a monoid by the vertices of $P$
\end{enumerate}   
then the conclusion of Theorem \ref{cor:main} holds where, in the 
statement of the theorem, $n$ has the value $L(\beta)$ for the unique 
minimal element $\beta$ of $\bar{E}_\Phi$.
\label{rem:compressed}
\end{corollary}
\begin{proof}
Let $\beta$ be the unique minimal element of $\bar{E}_\Phi$, whose 
existence is guaranteed by (ii). Then (iii) implies that $\beta = v_1 
+ v_2 + \cdots + v_n$ for some vertices $v_1$, $v_2,\dots,v_n$ of $P$, 
which must be pairwise distinct. Because of (i) any ordering $(v_p, 
v_{p-1},\dots,v_1)$ of the vertices of $P$ satisfies the assumptions
of Corollary \ref{thm:main}. The result follows from this corollary
observing that $L(\beta) = n$.
\end{proof}

General conditions on $\Phi$ and $P$ which guarantee assumptions (i) 
and (iii) of Corollary \ref{rem:compressed} were given by Ohsugi and 
Hibi \cite{OH}.
\begin{corollary}
Let $P$ be an $m$-dimensional polytope, as in {\rm (\ref{P})}. If $P$ 
is a {\rm 0-1} polytope and $R_\Phi$ is Gorenstein then the conclusion 
of Corollary \ref{rem:compressed} holds.
\label{cor:OH1}
\end{corollary}
\begin{proof}
Since $P$ is a 0-1 polytope, Theorem 1.1 and Lemma 2.1 in \cite{OH} imply,
respectively, that $P$ is compressed and that $E_\Phi$ is generated as 
a monoid by the integer points of $P$, which are exactly the vertices of
$P$. Thus the assumptions of Corollary \ref{rem:compressed} hold.
\end{proof}

\begin{remark} {\rm
Recall that a matrix is called \emph{totally unimodular} if all its 
subdeterminants are equal to $0, -1$ or $1$. Suppose that the linear 
functional $L(x) = a_1 x_1 + a_2 x_2 + \cdots + a_q x_q$ satisfies $a_i 
\in \ZZ$ for all $i$. It follows from the results of \cite{HK} that 
$P$ is an integer polytope if the matrix obtained from $\Phi$ by adding 
the row $(a_1, a_2,\dots,a_q)$ is totally unimodular. Hence, in view of
Corollary \ref{cor:OH1}, this statement and the assumptions that (i) 
$P \subseteq [0, 1]^q$ and (ii) $R_\Phi$ is Gorenstein imply the 
conclusion of Corollary \ref{thm:main}.
} 
\label{cor:OH2}
\end{remark}

\vspace{0.1 in}
\noindent
{\bf Magic labelings of graphs.} Let $G$ be a graph (multiple edges and 
loops allowed) with $p$ vertices and $q$ edges and edge set $\eE$. A 
\emph{magic labeling} \cite{St1} of $G$ of \emph{index} $r$ is an 
assignment $\ell: \eE \to \NN$ of nonnegative integers to the edges of 
$G$ such that for each vertex $v$ of $G$ the sum of the labels of all 
edges incident to $v$ is equal to $r$, in other words,
\[ \sum_{e: \ v \in e} \ \ell(e)  = r. \]
Let $\RR^\eE$ denote the real vector space with basis $\eE$ and let
$x_e$ be the linear functional on $\RR^\eE$ dual to the basis element 
$e \in \eE$. Let $\Phi$ be a full rank, $p' \times q$ integer matrix 
with kernel the set of $x \in \RR^\eE$ satisfying the linear system 
of equations of the form
\begin{equation}
\sum_{e: \ v \in e} \ x_e  = \sum_{e': \ v' \in e'} \ x_{e'}
\label{inc}
\end{equation}
in $\RR^\eE$, where $v, v'$ are vertices of $G$. If $L(x)$ is any of 
the functionals on $\RR^\eE$ in (\ref{inc}) and $P$ is as in (\ref{P}) 
then $\erh(P, r)$ counts the number $H_G (r)$ of magic labelings of 
$G$ of index $r$. It follows from \cite[Proposition 2.9]{St1} and either 
\cite[Theorem 2.3]{St3} (applied as in \cite[Example 2.4 (b)]{St3} in the 
case of the Birkhoff polytope) or \cite[Theorem 1.1]{OH} that conditions 
(i) and (iii) of Corollary \ref{rem:compressed} are both satisfied if 
the graph $G$ is bipartite (or, more generally, if it satisfies 
condition (iii) of \cite[Proposition 2.9]{St1}). Clearly condition (ii) 
is satisfied if $G$ is regular, since then $(1, 1,\dots,1) \in E_\Phi$. 
Assuming further that $G$ is connected we have 
$p' = p-2$. The following corollary specializes to Corollary 
\ref{cor:magic} when $G$ is the complete bipartite graph on two sets 
of vertices, each of size $n$.  
\begin{corollary}
For $n \ge 1$ and for any connected regular bipartite graph $G$ with 
$p$ vertices and $q=np/2$ edges we have
\[ \sum_{r \ge 0} \, H_G (r) \, t^r = \frac{h(t)} {(1 - t)^{m + 1}} \]
where $m = q - p + 1$ and $h(t) = h_0 + h_1 t + \cdots + h_d t^d$ is the 
$h$-polynomial of the boundary complex of a simplicial polytope of dimension 
$d = m-n+1$, so that $h(t)$ satisfies the conditions in the $g$-theorem. 

In particular  $h_i = h_{d-i}$ for all $i$ and $1 = h_0 \le h_1 \le \cdots 
\le h_{\lfloor d/2 \rfloor}$.
\label{cor:graphs}
\end{corollary}

\end{document}